\documentclass[11pt,reqno]{amsart}
\usepackage{enumerate, latexsym, amsmath, amsfonts, amssymb, amsthm, color}
\def\pmod #1{\ ({\rm{mod}}\ #1)}
\def\Z{\Bbb Z}

\def\l{\left}
\def\r{\right}
\def\bg{\bigg}
\def\({\bg(}
\def\){\bg)}
\def\t{\text}
\def\f{\frac}
\def\sign{\mathrm{sign}}

\def\per{{\rm per}}

\def\ls{\leqslant}

\def\bi{\binom}

\def\eq{\equiv}

\theoremstyle{plain}
\newtheorem{theorem}{Theorem}

\newtheorem{corollary}[theorem]{Corollary}
\newtheorem{conjecture}[theorem]{Conjecture}

\theoremstyle{remark}
\newtheorem{remark}[theorem]{Remark}

\allowdisplaybreaks

 \vspace{4mm}

\begin{document}

\hbox{Acta Math. Sinica Chin. Ser. 67 (2024), no.\,2, 286--295.}
\medskip

\title
[{On some determinants and permanents}]
{On some determinants and permanents}

\author
[Zhi-Wei Sun] {Zhi-Wei Sun}

\dedicatory{In memory of Prof. Yuan Wang}

\address{Department of Mathematics, Nanjing
University, Nanjing 210093, People's Republic of China}
\email{zwsun@nju.edu.cn}

\keywords{Determinant, permanent, $p$-adic congruence, Legendre symbol.
\newline \indent 2020 {\it Mathematics Subject Classification}. Primary 11C20, 15A15; Secondary 11A07, 11A15.
\newline \indent The work is supported
by the National Natural Science Foundation of China (grant 12371004).}

\begin{abstract}
In this paper we study some determinants and permanents. In particular, we investigate the new type determinants
$$\det[(i^2+cij+dj^2)^{p-2}]_{1\le i,j\le p-1}\  \text{and} \ \det[(i^2+cij+dj^2)^{p-2}]_{0\le i,j\le p-1}$$
modulo an odd prime $p$, where $c$ and $d$ are integers. We also pose some conjectures for further research.
\end{abstract}
\maketitle

\section{Introduction}
\setcounter{equation}{0}
\setcounter{theorem}{0}

Let $A=[a_{i,j}]_{1\ls i,j\ls n}$
be an matrix over a commutative ring. The determinant of $A$ and the permanent of $A$ are given by
$$\det(A)=\sum_{\sigma\in S_n}\sign(\sigma)\prod_{i=1}^na_{i,\sigma(i)}
\ \ \t{and}\ \ \per(A)=\sum_{\sigma\in S_n}\prod_{i=1}^na_{i,\sigma(i)}$$
respectively, where $S_n$ is the symmetric group of all permutations of $\{1,\ldots,n\}$,
and $\sign(\sigma)$ is $1$ or $-1$ according as the permutation $\sigma\in S_n$ is even or odd.

For $i,j=1,\ldots,n$, define
$$p_{ij}=\begin{cases}1&\t{if}\ i+j\ \t{is prime},
\\0&\t{otherwise}.\end{cases}$$
It is known that $|\det[p_{ij}]_{1\ls i,j\ls n}|$ is always a perfect square
(cf. \cite{Sp}). Our following theorem is a further extension of this result.

\begin{theorem}\label{Th1.3} Let $A=[a_{i,j}]_{1\ls i,j\ls n}$
be an matrix over a commutative ring. Suppose that $a_{i,j}=0$ whenever $i+j$ is an even number greater than two.

{\rm (i)} If $n=2m$ for some $m\in\Z^+$, then
\begin{equation}\label{2m-per}\per(A)=\per[a_{2i,2j-1}]_{1\ls i,j\ls m}\,\per[a_{2i-1,2j}]_{1\ls i,j\ls m}
\end{equation}
and
\begin{equation}\label{2m}\det(A)=(-1)^m\det[a_{2i,2j-1}]_{1\ls i,j\ls m}\,\det[a_{2i-1,2j}]_{1\ls i,j\ls m}.
\end{equation}

{\rm (ii)} If $n=2m+1$ for some $m\in\Z^+$, then
\begin{equation}\label{2m+1-per}\per(A)=a_{1,1}\,\per[a_{2i,2j+1}]_{1\ls i,j\ls m}\,\per[a_{2i+1,2j}]_{1\ls i,j\ls m}
\end{equation}
and
\begin{equation}\label{2m+1}\det(A)=(-1)^ma_{1,1}\det[a_{2i,2j+1}]_{1\ls i,j\ls m}\,\det[a_{2i+1,2j}]_{1\ls i,j\ls m}.
\end{equation}
\end{theorem}

Clearly, Theorem \ref{Th1.3} has the following consequence.

\begin{corollary} Let $A=[a_{i,j}]_{1\ls i,j\ls n}$
be a symmetric matrix over a commutative ring, with $a_{i,j}=0$
whenever $i+j\in\{2q:\ q=2,3,\ldots\}$.

{\rm (i)} If $n$ is even, then
$$\per(A)=\per[a_{2i,2j-1}]_{1\ls i,j\ls n/2}^2$$
and
$$ (-1)^{n/2}\det(A)=\det[a_{2i,2j-1}]_{1\ls i,j\ls n/2}^2.$$

{\rm (ii)} If $n>1$ and $2\nmid n$, then
$$\per(A)=a_{11}\,\per[a_{2i,2j+1}]_{1\ls i,j\ls (n-1)/2}^2$$
and
$$ (-1)^{(n-1)/2}\det(A)=a_{1,1}\det[a_{2i,2j-1}]_{1\ls i,j\ls (n-1)/2}^2.$$
\end{corollary}
\begin{remark} By B. Cloitre's comments linked in \cite{Sp}, one of his colleagues proved in 2002 that
if $A=[a_{ij}]_{1\ls i,j\ls n}$ is a symmetric matrix over a commutative ring with $a_{11}$ a square
and $a_{ij}=0$ whenever $i+j\in\{2q:\ q=2,3,\ldots\}$, then $|\det(A)|$ is a square.
\end{remark}

A matrix $A=[a_{ij}]_{1\ls i,j\ls n}$ over a commutative ring is called {\it skew-symmetric}
if $a_{ji}=-a_{ij}$ for all $i,j=1,\ldots,n$. For such a matrix $A$, clearly $2\det(A)=2\per(A)=0$
if $n$ is odd.

Theorem \ref{Th1.3} also has the following consequence.

\begin{corollary} Let $A=[a_{ij}]_{1\ls i,j\ls 2m}$ be a skew-symmetric matrix of even order
over a commutative ring such that $a_{ij}=0$ whenever $i+j\in\{2q:\ q=2,3,\ldots\}$.
Then
$$\per(A)=(-1)^m\per[a_{2i,2j-1}]_{1\ls i,j\ls m}^2$$
and
$$\det(A)=\det[a_{2i,2j-1}]_{1\ls i,j\ls m}^2.$$
\end{corollary}
\begin{remark} By a theorem of Cayley, any skew-symmetric determinant
over $\Z$ of even order is a perfect square (cf. \cite[[pp.\,23-24]{K}).
\end{remark}

Let $p$ be an odd prime and let $(\f{\cdot}p)$ be the Legendre symbol.
If $p\eq 3\mod 4$, then $i^2+j^2\not\eq0\pmod p$ for all $i,j\in\Z$, and
the author \cite[Theorem 1.4(ii)]{S19} proved that
$$\det\l[\f1{i^2+j^2}\r]_{1\ls i,j\ls(p-1)/2}\eq\l(\f 2p\r)\pmod p.$$
In the case $p\eq2\pmod3$, the author \cite[Remark 1.3]{S19} conjectured that
$$\det\l[\f1{i^2-ij+j^2}\r]_{1\ls i,j\ls p-1}\eq2x^2\pmod p$$
for some integer $x\not\eq0\pmod p$.
This was recently confirmed by Wu, She and Ni \cite{WSN}.
Note that for any integer $x\not\eq0\pmod p$ we have $\f1x\eq x^{p-2}\pmod p$
by Fermat's little theorem.

Let $P(x,y)\in F[x,y]$ with $F$ a field. If
$$P(x,j)=\sum_{k=0}^{n-1}a_{jk}x^k\quad\t{for all}\ j=1,\ldots,n,$$
then
$$\det[P(x_i,j)]_{1\ls i,j\ls n}=\det[a_{jk}]_{1\ls j\ls n\atop 0\ls k<n}\times\prod_{1\ls i<j\ls n}(x_j-x_i)$$
by a known result (cf. \cite{K}). It follows that
$$\deg P<n-1\ \Longrightarrow\ \det[P(i,j)]_{1\ls i,j\ls n}=0.$$

Now we state our second theorem.

\begin{theorem}\label{Th1.1}
Let $c,d\in\Z$. For any prime $p>3$, we have
\begin{equation}
\label{0-p-1}\det[(i^2+cij+dj^2)^{p-2}]_{0\ls i,j\ls p-1}\eq0\pmod p.
\end{equation}
\end{theorem}
\begin{remark} For $p=3$ and $c,d\in\Z$, it is easy to find that
$$\det[(i^2+cij+dj^2)^{p-2}]_{0\ls i,j\ls p-1}=\det[i^2+cij+dj^2]_{0\ls i,j\ls 2}=-4cd.$$
\end{remark}

\noindent {\bf Definition 1.1}. For any $c,d\in\Z$ and odd prime $p$, we set
\begin{equation}\label{D_p(c,d)}
D_p(c,d):=\det[(i^2+cij+dj^2)^{p-2}]_{1\ls i,j\ls p-1}.
\end{equation}

Let $c,d\in\Z$ and let $p$ be an odd prime. As
\begin{align*}&\det[((p-i)^2+c(p-i)j+dj^2)^{p-2}]_{1\ls i,j\ls p-1}
\\=&(-1)^{\sum_{k=1}^{p-2}k}\det[(i^2+cij+dj^2)^{p-2}]_{1\ls i,j\ls p-1},
\end{align*}
we have
\begin{equation}
\label{-c,d}D_p(-c,d)\eq\l(\f{-1}p\r)D_p(c,d)\pmod p.
\end{equation}

\begin{theorem}\label{Th1.2} Let $p>3$ be a prime.

{\rm (i)} If $p\eq3\pmod4$, then
$$D_p(c,-1)\eq0\pmod p\quad\t{for all}\ c\in\Z.$$

{\rm (ii)} When $p\eq3\pmod4$, we have
$$D_p(2,2)\eq0\pmod p.$$
If $p\eq\pm1\pmod{12}$, then
$$D_p(6,6)\eq0\pmod p.$$
\end{theorem}

Theorem 1.1, and Theorems 1.6 and 1.8, will be proved in Sections 2 and 3, respectively.
In Section 4, we pose some conjectures on determinants and permanents.

\section{Proofs of Theorem 1.1}
\setcounter{equation}{0}
\setcounter{theorem}{0}

\medskip
\noindent{\tt Proof of Theorem 1.1(i)}. Suppose that $n=2m$ with $m\in\Z^+$.
Note that
$$\per(A)=\sum_{\sigma\in S_{2m}}\prod_{i=1}^{2m}a_{i,\sigma(i)}
\ \t{and}\ \det(A)=\sum_{\sigma\in S_{2m}}\sign(\sigma)\prod_{i=1}^{2m}a_{i,\sigma(i)}.$$

Assume that $\sigma\in S_{2m}$ and $\prod_{i=1}^{2m}a_{i,\sigma(i)}\not=0$. Then there are $\tau_1,\tau_2\in S_m$ such that for any $i,j=1,\ldots,m$ we have
$$\sigma(2i)=2\tau_1(i)-1\ \ \t{and}\ \ \sigma(2j-1)=2\tau_2(j).$$
For $i=1,\ldots,m$ define $\sigma'(i)=\sigma(2i)$ and $\sigma'(m+i)=\sigma(2i-1)$. It is easy to see that
$\sigma'=\rho\sigma$ for some $\rho\in S_{2m}$ which is a product of $1+2+\cdots+m$ transposes. Thus
$$\sign(\sigma')=(-1)^{m(m+1)/2}\sign(\sigma).$$
Note that $\sigma'(i)=2\tau_1(i)-1$ and $\sigma'(m+i)=2\tau_2(i)$ for all $i=1,\ldots,m$.
To find the sign of $\sigma'$, we investigate the parity of its inverse pairs.
For $1\ls i<j\ls m$, clearly
 $$\sigma'(i)>\sigma'(j)\iff 2\tau_1(i)-1>2\tau_1(j)-1\iff \tau_1(i)>\tau_1(j),$$
and
 $$\sigma'(m+i)>\sigma'(m+j)\iff 2\tau_2(i)>2\tau_2(j)\iff \tau_2(i)>\tau_2(j).$$
Moreover,
 \begin{align*}&|\{(i,m+j):\ 1\ls i,j\ls m\ \t{and}\ \sigma'(i)>\sigma'(m+j)\}|
 \\=\ &|\{(i,j):\ 1\ls i,j\ls m\ \t{and}\ 2\tau_1(i)-1>2\tau_2(j)\}|
 \\=\ &|\{(\tau_1(i),\tau_2(j)):\ 1\ls i,j\ls m\ \t{and}\ 2\tau_1(i)-1>2\tau_2(j)\}|
 \\=\ &|\{(s,t):\ 1\ls s,t\ls m\ \t{and}\ 2s-1>2t\}|=\sum_{s=1}^m\sum_{0<t<s}1=\sum_{s=1}^m(s-1).
 \end{align*}
 Thus,
 $$\sign(\sigma')=(-1)^{\sum_{s=1}^m(s-1)}\sign(\tau_1)\sign(\tau_2)=(-1)^{m(m+1)/2-m}\sign(\tau_1)\sign(\tau_2),$$
 and hence
 $$\sign(\sigma)=(-1)^{m(m+1)/2}\sign(\sigma')=(-1)^m\sign(\tau_1)\sign(\tau_2).$$
Therefore
\begin{align*}
 \prod_{i=1}^{2m}a_{i,\sigma(i)}=&\prod_{i=1}^ma_{2i,2\tau_1(i)-1}\times\prod_{j=1}^ma_{2j-1,2\tau_2(j)}
 \end{align*}
 and
 \begin{align*}
 \sign(\sigma)\prod_{i=1}^{2m}a_{i,\sigma(i)}=&(-1)^m\sign(\tau_1)\prod_{i=1}^ma_{2i,2\tau_1(i)-1}\times
 \sign(\tau_2)\prod_{j=1}^ma_{2j-1,2\tau_2(j)}.
\end{align*}

 In view of the above, we have
  \begin{align*}\per(A)=&\sum_{\sigma\in S_{2m}}\prod_{i=1}^{2m}a_{i,\sigma(i)}
  \\=&\(\sum_{\tau_1\in S_m}\prod_{i=1}^ma_{2i,2\tau_1(i)-1}\)\sum_{\tau_2\in S_m}\prod_{j=1}^m a_{2j-1,2\tau_2(j)}
 \\=&\per[a_{2i,2j-1}]_{1\ls i,j\ls m}\times\per[a_{2i-1,2j}]_{1\ls i,j\ls m}
  \end{align*}
  and
  \begin{align*}\det(A)=&\sum_{\sigma\in S_{2m}}\sign(\sigma)\prod_{i=1}^{2m}a_{i,\sigma(i)}
  \\=&(-1)^m\sum_{\tau_1\in S_m}\sign(\tau_1)a_{2i,2\tau_1(i)-1}\sum_{\tau_2\in S_m}\sign(\tau_2)\prod_{j=1}^m a_{2j-1,2\tau_2(j)}
  \\=&(-1)^m\det[a_{2i,2j-1}]_{1\ls i,j\ls m}\times\det[a_{2i-1,2j}]_{1\ls i,j\ls m}.
  \end{align*}
  This concludes the proof. \qed

\medskip
\noindent{\tt Proof of Theorem 1.1(ii)}. Suppose that $n=2m+1$ with $m\in\Z^+$.
Then
$$\per(A)=\sum_{\sigma\in S_{2m+1}}\prod_{i=1}^{2m+1}a_{i,\sigma(i)}
\ \t{and}\ \det(A)=\sum_{\sigma\in S_{2m+1}}\sign(\sigma)\prod_{i=1}^{2m+1}a_{i,\sigma(i)}.$$

Let $\sigma\in S_{2m+1}$ with $\prod_{i=1}^{2m+1}a_{i,\sigma(i)}\not=0$. Then
$$\{i+\sigma(i):\ i=1,\ldots,2m+1\}\cap\{2q:\ q=2,\ldots,m\}=\emptyset.$$
If $\sigma(1)\not=1$, then those $\sigma(2i+1)\ (0\ls i\ls m)$ are pairwise distinct even numbers,
which is impossible since $|\{2j:\ j=1,\ldots,m\}|<m+1$. Thus $\sigma(1)=1$.
As $\sigma(2i)\ (1\ls i\ls m)$ are pairwise distinct odd numbers, for some $\tau_1\in S_m$
we have $\sigma(2i)=2\tau_1(i)+1$ for all $i=1,\ldots,m$. As $\sigma(2i+1)\ (1\ls i\ls m)$ are distinct
positive even integers, for some $\tau_2\in S_m$ we have $\sigma(2i+1)=2\tau_2(i)$
for all $i=1,\ldots,m$. Let $\sigma'(1)=1$, $\sigma'(i+1)=\sigma(2i)$
and $\sigma'(m+1+i)=\sigma(2i+1)$ for all $i=1,\ldots,m$.
It is easy to see that $\sigma'=\lambda\sigma$ for some $\lambda\in S_{2m+1}$ which is a product of
$\sum_{0<k<m}k=m(m-1)/2$ transposes. Thus $$\sign(\sigma')=(-1)^{m(m-1)/2}\sign(\sigma).$$
Note that $\sigma'(i+1)=2\tau_1(i)+1$ and $\sigma'(m+1+i)=2\tau_2(i)$ for all $i=1,\ldots,m$.
When $1\ls i<j\ls m$, we have
 $$\sigma'(i+1)>\sigma'(j+1)\iff 2\tau_1(i)+1>2\tau_1(j)+1\iff \tau_1(i)>\tau_1(j),$$
and
 $$\sigma'(m+1+i)>\sigma'(m+1+j)\iff 2\tau_2(i)>2\tau_2(j)\iff \tau_2(i)>\tau_2(j).$$
Moreover,
 \begin{align*}&|\{(i+1,m+1+j):\ 1\ls i,j\ls m\ \t{and}\ \sigma'(i+1)>\sigma'(m+1+j)\}|
 \\=\ &|\{(i,j):\ 1\ls i,j\ls m\ \t{and}\ 2\tau_1(i)+1>2\tau_2(j)\}|
 \\=\ &|\{(\tau_1(i),\tau_2(j)):\ 1\ls i,j\ls m\ \t{and}\ 2\tau_1(i)+1>2\tau_2(j)\}|
 \\=\ &|\{(s,t):\ 1\ls s,t\ls m\ \t{and}\ 2s+1>2t\}|=\sum_{s=1}^m\sum_{t=1}^s1=\sum_{s=1}^m s.
 \end{align*}
Thus
$$\sign(\sigma')=(-1)^{\sum_{s=1}^ms}\sign(\tau_1)\sign(\tau_2)=(-1)^{m(m-1)/2+m}\sign(\tau_1)\sign(\tau_2),$$
and hence
$$\sign(\sigma)=(-1)^{m(m-1)/2}\sign(\sigma')=(-1)^m\sign(\tau_1)\sign(\tau_2).$$
Note that
$$\prod_{i=1}^{2m+1}a_{i,\sigma(i)}
=a_{1,1}\prod_{i=1}^ma_{2i,2\tau_1(i)+1}\times\prod_{i=1}^ma_{2i+1,2\tau_2(i)}
$$
and
\begin{align*}&\sign(\sigma)\prod_{i=1}^{2m+1}a_{i,\sigma(i)}
=(-1)^ma_{1,1}\sign(\tau_1)\prod_{i=1}^ma_{2i,2\tau_1(i)+1}\times\sign(\tau_2)\prod_{i=1}^ma_{2i+1,2\tau_2(i)}.
\end{align*}

In light of the above, we see that
\begin{align*}\per(A)=&a_{1,1}\(\sum_{\tau_1\in S_m}\prod_{i=1}^m a_{2i,2\tau_1(i)+1}\)
\sum_{\tau_2\in S_m}\prod_{i=1}^ma_{2i+1,2\tau_2(i)}
\\=&a_{1,1}\per[a_{2i,2j+1}]_{1\ls i,j\ls m}\times\per[a_{2i+1,2j}]_{1\ls i,j\ls m}
\end{align*}
and
\begin{align*}\det(A)=&(-1)^ma_{1,1}\(\sum_{\tau_1\in S_m}\sign(\tau_1)\prod_{i=1}^m a_{2i,2\tau_1(i)+1}\)
\sum_{\tau_2\in S_m}\sign(\tau_2)\prod_{i=1}^ma_{2i+1,2\tau_2(i)}
\\=&(-1)^ma_{1,1}\det[a_{2i,2j+1}]_{1\ls i,j\ls m}\times\det[a_{2i+1,2j}]_{1\ls i,j\ls m}.
\end{align*}
This concludes the proof. \qed

\section{Proofs of Theorems 1.6 and 1.8}
\setcounter{equation}{0}
\setcounter{theorem}{0}

\medskip
\noindent{\tt Proof of Theorem 1.6}.
Let $a_{ij}=(i^2+cij+dj^2)^{p-2}$ for all $i,j=0,\ldots,p-1$.
If $p\mid d$, then $a_{0j}\eq0\pmod p$ for all $j=0,\ldots,p-1$, and hence $\det[a_{ij}]_{0\ls i,j\ls p-1}\eq0\pmod p$.

Below we assume that $p\nmid d$.
Fix $j\in\{1,\ldots,p-1\}$. Then
\begin{align*}4^{p-2}\sum_{i=0}^{p-1}a_{ij}=&\sum_{i=0}^{p-1}(4i^2+4cij+4dj^2)^{p-2}
=\sum_{i=0}^{p-1}((2i+cj)^2+(4d-c^2)j^2)^{p-2}
\\\eq&\sum_{i=1}^{p}(i^2+(4d-c^2)j^2)^{p-2}=\sum_{i=1}^{p}\sum_{k=0}^{p-2}\bi{p-2}ki^{2k}((4d-c^2)j^2)^{p-2-k}
\\=&\sum_{k=0}^{p-2}\bi{p-1}{k+1}\f{k+1}{p-1}(4d-c^2)^{p-2-k}j^{2(p-2-k)}\sum_{i=1}^{p}i^{2k}
\pmod p.
\end{align*}
Clearly $\sum_{i=1}^pi^0=p\eq0\pmod p$. For $k\in\{1,\ldots,p-1\}$ with $k\not=(p-1)/2$, as $p-1\nmid2k$ we have
$$\sum_{i=1}^pi^{2k}\eq\sum_{i=1}^{p-1}i^{2k}\eq0\pmod p$$
(cf. \cite{IR}). For $k=(p-1)/2$, we have
$$\sum_{i=1}^pi^{2k}\eq\sum_{i=1}^{p-1}i^{p-1}\eq\sum_{i=1}^{p-1}1\eq-1\pmod p$$
by Fermat's little theorem.
Therefore,
\begin{align*}4^{p-2}\sum_{i=0}^{p-1}a_{ij}\eq&-\bi{p-1}{(p+1)/2}\f{(p+1)/2}{p-1}(4d-c^2)^{p-2-(p-1)/2}j^{2(p-2-(p-1)/2)}
\\\eq&\f12(-1)^{(p+1)/2}(4d-c^2)^{(p-1)/2-1}j^{p-3}\pmod p
\end{align*}
and hence
$$\sum_{i=0}^{p-1}a_{ij}\eq\f{2}{j^2}(c^2-4d)^{(p-3)/2}\pmod p$$
by Fermat's little theorem. As $a_{0j}=(dj^2)^{p-2}\eq1/(dj^2)\pmod p$, we have
\begin{equation}\label{2.1}\l(1-2d(c^2-4d)^{(p-3)/2}\r)a_{0j}+\sum_{i=1}^{p-1}a_{ij}\eq0\pmod p.
\end{equation}

Note that
$$\l(1-2d(c^2-4d)^{(p-3)/2}\r)a_{00}+\sum_{i=1}^{p-1}a_{i0}
=\sum_{i=1}^{p-1}i^{2(p-2)}\eq\sum_{i=1}^{p-1}\f1{i^2}\eq0\pmod p$$
by Wolstenholme's theorem. Combining this with the above paragraph, we see that
$$\l(1-2d(c^2-4d)^{(p-3)/2}\r)a_{0j}+\sum_{i=1}^{p-1}a_{ij}\eq0\pmod p
\quad\t{for all}\ j=0,\ldots,p-1.$$
It follows that $$\det[a_{ij}]_{0\ls i,j\ls p-1}\eq0\pmod p$$
as desired.

By the above, we have proved Theorem 1.6. \qed

\medskip
\noindent
{\tt Proof of Theorem 1.8}. (i) Assume that $p\eq3\pmod4$. For any $c\in\Z$, we have
\begin{align*}D_p(c,-1)=&\det[(j^2+cij-i^2)^{p-2}]_{1\ls i,j\ls p-1}
\\=&\det[-(i^2-cij-j^2)^{p-2}]_{1\ls i,j\ls p-1}=D_p(-c,-1)
\\\eq&\l(\f{-1}p\r)D_p(c,-1)=-D_p(c,-1)\pmod p
\end{align*}
with the help of \eqref{-c,d}, therefore $D_p(c,-1)\eq0\pmod p$.

(ii) Let $c,d\in\Z$ with $p\nmid d(c^2-4d)$. In view of \eqref{2.1}, for any $j=1,\ldots,p-1$ we have
\begin{align*}\sum_{i=1}^{p-1}(i^2+cij+dj^2)^{p-2}
\eq&\l(2d(c^2-4d)^{(p-3)/2}-1\r)\f1{dj^2}
\\\eq&\l(\f{2d}{c^2-4d}\l(\f{c^2-4d}p\r)-1\r)\f1{dj^2}
\pmod p.\end{align*}
So, $D_p(c,d)\eq0\pmod p$ provided that
$$\f{c^2-4d}{2d}\eq\l(\f{c^2-4d}p\r)\pmod p.$$
Hence $D_p(2,2)\eq0\pmod p$
if $p\eq3\pmod4$, and $D_p(6,6)\eq0\pmod p$
if $p\eq\pm1\pmod{12}.$

In view of the above, we have completed the proof of Theorem 1.8. \qed

\section{Some conjectures}
\setcounter{equation}{0}
\setcounter{theorem}{0}

\begin{conjecture} Let $n>3$ be an odd integer, and let $c,d\in\Z$ with the Jacobi symbol $(\f dn)$ equal to $-1$. Then
$$\det[(i^2+cij+dj^2)^{n-2}]_{0\ls i,j\ls n-1}\eq0\pmod{n^2}.$$
\end{conjecture}

\begin{conjecture} For any prime $p\eq1\pmod4$ with $p\eq\pm2\pmod5$,  we have
$$\l(\f{D_p(1,-1)}p\r)=1.$$
\end{conjecture}

\begin{conjecture} For any odd prime $p$, we have
$$\l(\f{D_p(2,-1)}p\r)=-1\iff p\eq5\pmod8.$$
\end{conjecture}

\begin{conjecture} For any odd prime $p\eq\pm2\pmod5$, we have
$$\l(\f{D_p(3,1)}p\r)=\begin{cases}(\f 6p)&\t{if}\ p\eq1\pmod4,
\\0&\t{if}\ p\eq3\pmod4.\end{cases}$$
\end{conjecture}

\begin{conjecture} Let $p>3$ be a prime with $(\f p7)=-1$ and $p\not\eq15\pmod{16}$. Then
$$\l(\f{D_p(1,16)}p\r)=\l(\f{-2}p\r).$$
\end{conjecture}

\begin{conjecture} Let $p$ be a prime. If $p\eq5\pmod{24}$, then $p\nmid D_p(6,6)$.
 If $p\eq 19\pmod{24}$, then $(\f{D_p(6,6)}p)\not=-1$.
\end{conjecture}

For any positive integer $n$, we let $D(n)$ denote the set of all derangements of $1,\ldots,n$, i.e.,
$$D(n)=\{\tau\in S_n:\ \tau(j)\not=j\ \t{for all}\ j=1,\ldots,n\}.$$

\begin{conjecture} For any odd prime $p$, we have
$$\sum_{\tau\in D(p-1)}\prod_{j=1}^{p-1}\f1{j-\tau(j)}\eq\l(\f{-1}p\r)\pmod{p^2}$$
and
$$\sum_{\tau\in D(p-1)}\sign(\tau)\prod_{j=1}^{p-1}\f1{j-\tau(j)}\eq1\pmod{p^2}.$$
\end{conjecture}

\begin{conjecture} {\rm (i)} For any odd prime $p$ we have
$$\sum_{\tau\in D(p-1)}\prod_{j=1}^{p-1}\f{j+\tau(j)}{j-\tau(j)}\eq1-2\l(\f{-1}p\r)\pmod{p}.$$

{\rm (ii)} For any prime $p>3$, the number
$$\f1{p^{3-(\f{-1}p)}}\sum_{\tau\in D(p-1)}\sign(\tau)\prod_{j=1}^{p-1}\f{j+\tau(j)}{j-\tau(j)}$$
is a quadratic residue modulo $p$.
\end{conjecture}

For convenience, we view an empty product $\prod_{i\in\emptyset}a_i$ as $1$.

\begin{conjecture} Let $p$ be an odd prime. Then
$$\sum_{\tau\in S_{p-1}}\prod_{j=1\atop \tau(j)\not=j}^{p-1}\frac1{j-\tau(j)}\equiv 1+\l(\f{-1}p\r)\pmod p.$$
Moreover, when $p\eq3\pmod4$ we have
$$\sum_{\tau\in S_{(p-1)/2}}\prod_{j=1\atop \tau(j)\not=j}^{(p-1)/2}\frac1{j^2-\tau(j)^2}
\equiv 1\pmod p.$$
\end{conjecture}

\begin{conjecture} Let $p$ be an odd prime. Then
$$\sum_{\tau\in S_{p}}\prod_{j=1\atop \tau(j)\not=j}^{p}\frac{j+\tau(j)}{j-\tau(j)}\equiv1-\l(\f{-1}p\r)\pmod p$$
and
$$\sum_{\tau\in S_{p}}\sign(\tau)\prod_{j=1\atop \tau(j)\not=j}^{p}\frac{j+\tau(j)}{j-\tau(j)}\equiv-\f p2\pmod {p^2}.$$
\end{conjecture}

\begin{conjecture} {\rm (i)} For any prime $p$, we have
$$\sum_{\tau\in S_{p-1}}\prod_{j=1\atop \tau(j)\not=j}^{p-1}\frac{j+\tau(j)}{j-\tau(j)}\equiv((p-2)!!)^2\pmod{p^2}.$$

{\rm (ii)} If $p$ is an odd prime, then
$$\sum_{\tau\in S_{p-1}}\mathrm{sign}(\tau)\prod^{p-1}_{j=1\atop \tau(j)\not=j}\frac{j+\tau(j)}{j-\tau(j)}\equiv\frac{(-1)^{(p+1)/2}}{p-2}((p-2)!!)^2\pmod{p^2}.$$
\end{conjecture}

\begin{conjecture} Let $p>3$ be a prime. If $p\equiv3\pmod4$, then
$$\sum_{\tau\in S_{(p-1)/2}}\mathrm{sign}(\tau)\prod^{(p-1)/2}_{j=1\atop \tau(j)\not=j}\frac{j^2+\tau(j)^2}{j^2-\tau(j)^2}\equiv0\pmod{p^2}.$$
When $p\equiv 7\pmod 8$, we have
$$\sum_{\tau\in S_{(p-1)/2}}\mathrm{sign}(\tau)\prod^{(p-1)/2}_{j=1\atop \tau(j)\not=j}\frac{j^2+\tau(j)^2}{j^2-\tau(j)^2}\equiv0\pmod{p^3}.$$
\end{conjecture}

\end{document}